\documentclass{amsart}
\usepackage{amscd,amssymb,graphicx}
\author[Fialowski]{Alice Fialowski}
\address{
Alice Fialowski\\
E\"otv\"os Lor\'and University\\
Budapest, Hungary} \email{fialowsk@cs.elte.hu}
\author[Penkava]{Michael Penkava}
\address{
Michael Penkava\\
University of Wisconsin-Eau Claire\\
Eau Claire, WI 54702-4004} \email{penkavmr@uwec.edu}
\thanks{Research of the first author was partially supported by  OTKA grant K77757 and the second author by grants from the
University of Wisconsin-Eau Claire.}

\newtheorem{thm}{Theorem}[section]

\theoremstyle{definition}

\def \ph{\varphi}

\def\GL{\mathbb{GL}}

\def \refeq#1{equation (\ref{#1})}
\def \ra{\rightarrow}

\def \hom{\mbox{\rm Hom}}

\def \gl{\mbox{$\mathfrak{gl}$}}
\def \tns{\otimes}

\def \mcom{,\cdots,}

\def \C{\mbox{$\mathbb C$}}
\def \Z{\mbox{$\mathbb Z$}}

\def\zt{\mbox{$\Z_2$}}

\def\inv{^{-1}}

\def\im{\operatorname{Im}}

\def\m{\mbox{$\mathfrak m$}}

\def\V{V}

\def\and{\mbox{ \rm and }}

\DeclareMathOperator*{\invlim}{\underleftarrow{\rm lim}}

\def\inv{^{-1}}

\def\P{\mathbb P}

\input{xy}
\xyoption{all}

\begin{document}
\setlength{\multlinegap}{0pt}
\title[Versal deformations]
{On Arnold's versal deformations  and stratifications of moduli spaces given by group actions}

\address{}%
\email{}%

\thanks{}%

\date{\today}
\begin{abstract}
The aim of this paper is to compare stratifications of moduli spaces
given by group actions in the case of similarity of matrices introduced
by Arnold and the author's stratification by projective orbifolds, and
its relation to deformations of elements in the moduli space.
\end{abstract}
\maketitle

\section{Miniversal Deformations of algebras}
Let us consider the moduli space of associative or Lie algebras on a finite dimensional complex vector space $V$. First, the set of such algebra structures is a variety determined by a set of quadratic constraints on the coefficients of the algebra, expressed as an element of the vector space $C^2(V)$, where
$C^n(V)=\hom(V^{\tns n},V)$, in the case of associative algebras, and $C^n(V)=\hom(S^n(V),V)$ for Lie algebras. In either case, $C^n(V)$ is a finite dimensional complex vector space. Secondly, there is an action of the group $G=\GL(V)$ of linear automorphisms of $V$ on $C^2(V)$, which preserves the variety of algebras.  The moduli space of algebra structures on $V$ is precisely the set of equivalence classes of algebras under this group action.

In order to describe the deformation theory of algebras, we introduce the notion of a \emph{formal algebra}, that is, a complete local algebra, which is
a commutative, unital algebra with unique maximal ideal $\m$, satisfying $A=\invlim_n A/{\m}^n$, where equality means canonical isomorphism, such that $A/\m=\C$. An \emph{infinitesimal algebra} $A$ is a formal algebra for which ${\m}^2=0$. Given a formal algebra,  there is a unique morphism of algebras
$\epsilon:A\ra\C$, whose kernel is just $\m$.

A formal algebra $A$ is said to be the \emph{base} of a \emph{formal
deformation} $d_A$ of an algebra structure $d$, if we are given an
element $d_A\in C_A^2(V)=C^2(V)\tns A$, and $\epsilon_*(d_A)=d$, where $\epsilon_*(x\tns a)=\epsilon(a)x$, such that $d_A$ equips $V_A=V\tns A$ with the structure of
an algebra with coefficients in $A$. There is a natural decomposition $C_A^2(V)=C^2(V)\oplus C^2(V)\tns\m$ and in terms of this decomposition, a
deformation of $d$ with base $A$ is given by $d_A=d+\psi$ where $\psi\in C^2(V)\tns\m$.

Any $A$-linear map $f_A:V_A\ra V_A$ is determined by maps $f:V\ra V$ and $\ph:V\ra\V\tns\m$, so we can express $f_A=f+\ph$.

A formal automorphism $g_A$ of $V$ is an $A$-linear map $g_{A}:V_A\ra V_A$, such that $g_A=1_V+\ph$, where $\ph:V\ra V\tns\m$. A formal automorphism of $V$ is
automatically invertible.  Two formal deformations $d_A$ and $d'_A$ of $d$ are said to be \emph{formally equivalent} if there is a formal
automorphism $g_A$ of $V$ such that $g_A^*(d_A)=d'_A$, where $g_A^*(d_A)=g_A\inv d_A g_A$, and by the map $g_A$ on the right we mean the
map $C^2_A(V)\ra C^2_A(V)$ induced by $g_A$.

The space $C(V)=\prod_n C^n(V)$ is called the cochain complex of $V$ and it is equipped with the structure of a \zt-graded Lie algebra, which in the case of
associative algebras is called the Gerstenhaber bracket. We have
$[C^k(V),C^l(V)]\subseteq C^{k+l-1}(V)$, In terms of this bracket, the condition for $d\in C^(V)$ to be an algebra is simply
$[d,d]=0$, which is called the codifferential property. We define a differential $D$ on $C(V)$ by $D(\ph)=[d,\ph]$, and it is precisely the fact that
$d$ is a codifferential that guarantees that $D^2=0$. We have $D:C^k(V)\ra C^{k+1}(V)$ As a consequence, we can define the \emph{cohomology} of the
algebra by $$H^n(d)=
\ker(D:C^n(V)\ra C^{n+1}(V))/\im(D:C^{n-1}(V)\ra C^n(V).$$

The importance of this cohomology emerges when we consider \emph{infinitesimal} deformations of an algebra, which are formal deformations with an
infinitesimal base. An element $d_A=d+\psi$ where $\psi\in C^2(V)\tns\m$ is an infinitesimal deformation of $d$ precisely when the \emph{cocycle condition}
$D(\psi)=0$ is satisfied. This is because $[d_A,d_A]=[d,d]+2D(\psi)+[\psi,\psi]$, and we note that $[\psi,\psi]\in C^3(V)\tns\m^2$, so must vanish because
$\m^2=0$ for an infinitesimal algebra.

When $A$ is an infinitesimal algebra, a formal automorphism of $V$ is called an \emph{infinitesimal automorphism}, and a formal equivalence between two
infinitesimal deformations is called an \emph{infinitesimal equivalence}.

If $f:A\ra B$ is a morphism of two formal algebras, and $d_A$ is a deformation of $d$ with base $A$, then there is an induced deformation $f_*(d_A)$ with
base $B$ given as follows.  If $d_A=d+\ph$, then $d_B=d+f_*(\ph)$ where $f_*(\ph)=(1\tns f)\circ\ph$.  A deformation $d_A$ with base $A$ is said to be
\emph{versal} if given any deformation $d_B$ with base $B$, there is some $f:A\ra B$ such that $f_*(d_A)$ is formally equivalent to $d_B$. The deformation is
said to be \emph{universal}, if the map $f$ is unique. In general there is no universal deformation of $d$, but there is a type of deformation which has some
universal properties, called a \emph{miniversal} deformation.  A deformation $d_A$ is called miniversal, if it is versal and whenever $d_B$ is an
infinitesimal deformation, the morphism $f$ is unique. On algebraic
deformations see more in \cite{fi}, \cite{fi2}, \cite{fp10}.

The main theorem which relates cohomology to deformation theory is the existence of a universal deformation in the family of all infinitesimal deformations.
This type of universal deformation is guaranteed to exist whenever $H^2(d)$ is finite dimensional, which is always true if $V$ is a finite dimensional space.
\begin{thm}\cite{ff2}
Let $\langle \delta^1,\dots,\delta^m\rangle$ be a set of cocycles in $C^2(V)$, whose images give a basis of $H^2(d)$. Let $A=\C[t_1\mcom t_m]/(t_it_j)$ be
the infinitesimal base given by the quotient of the polynomial algebra on $m$ variables by the ideal $I$ generated by all quadratic polynomials $t_it_j$.
Then 
$$d_{\text{inf}}=d+t_i\delta^i$$
is a universal infinitesimal deformation.
\end{thm}
In this sense, the cohomology $H^2(d)$ classifies the infinitesimal deformations.  Moreover, one can show, \cite{ff2}, that a miniversal deformation 
of $d$ can be constructed with base given by a quotient of $\C[[t_1\mcom t_m]]$ by an ideal. The minimal number of parameters of a versal deformation of 
$d$ is precisely $\dim(H^2(d))$, and in \cite{ff2}, it was shown that there is a miniversal deformation with precisely this number of parameters.

\section{Arnold's presentation of versal deformations of matrices under a group action}
We want to present a picture of deformation theory in an special analytic context, given by an example presented by Arnold in \cite{arnold1},
where he considers deformations of the moduli space of similar matrices, that is $n\times n$ matrices under the action of $\GL(n)$ by conjugation.
It is well known that every matrix is conjugate to a matrix in Jordan normal form, and this form is unique up to the ordering of the Jordan blocks.

A deformation $A(\lambda)$ of an $n\times n$ matrix $A$ is a family of $n\times n$ matrices $A(\lambda)$ where the entries of the matrices are given by
power series in a family of parameters $\lambda_i$, which converge in a
neighborhood of the origin, such that $A(0)=A_0$. We can express this in a more algebraic language as follows. Let $B_k$ be the algebra
given by the germs of differentiable functions at the origin on $\C^k$ for some $k$. Then $B_k$ is a formal algebra, with maximal ideal $\m$ given by
the functions which vanish at the origin.
An element $A(\lambda)\in\gl(n)\tns B_k$ is a deformation of $A_0$ provided that $A(\lambda)=A_0+\Psi(\lambda)$, where $\Psi(\lambda)\in\gl(n)\tns\m$. If
$\ph$ is a germ of a holomorphic function $\C^l\ra\C^k$ at the origin, such that $\ph(0)=0$, and $A(\lambda)$ is a deformation of $A_0$ with base
$B_k$, then $\ph$ determines an induced deformation $\ph^*(A)(\mu)$ of $A_0$ with base
$B_l$, given by $\ph^*(A)(\mu)=A(\ph(\mu))$.

Two deformations $A(\lambda)$ and $B(\lambda)$ of $A_0$ with base $\lambda$ are equivalent if there is a deformation $C(\lambda)$ of the identity matrix,
such that $$B(\lambda)=C(\lambda)\inv A(\lambda)C(\lambda.$$ A deformation with base $B_k$ is said to be versal if given any deformation $Q(\mu)$ with
base $B_l$, there is a $\ph:C^l\ra C^k$ such that $Q(\mu)$ is equivalent to $\ph^*(A)(\mu)$. If $\ph$ is determined uniquely, then $A(\lambda)$ is called
a universal deformation of $A_0$.

Just as in the algebraic case, we cannot expect  that a universal deformation of $A_0$ generally exists. 
A versal deformation $A(\lambda)$ is called
miniversal when the dimension of the parameter space $B_k$ is minimal.  A versal deformation always exists, and therefore there is always
a miniversal deformation as well.

The notion of versality is similar but not identical to the algebraic case as the following example will show. Consider the matrix
$A_0=\left[\begin{smallmatrix}0&1\\0&0\end{smallmatrix}\right]$. Let $A(\lambda)$ be the deformation given by
$A(\lambda)=\left[\begin{smallmatrix}\lambda_1&1+\lambda_2\\0&+\lambda_3\end{smallmatrix}\right]$. This deformation was given in \cite{arnold1} as an example
of a deformation which is not versal.  To see this, suppose that $Q(\mu)=\left[\begin{smallmatrix}\mu_1&1+\mu_2\\\mu_3&\mu_4\end{smallmatrix}\right]$.
It is not hard to see that $Q(\mu)$ must be a versal deformation. Suppose we define $\ph:C^4\ra C^3$ by
$\ph(\mu_1,\mu_2,\mu_3)=(\lambda_1,\lambda_2,\lambda_3)$, where
\begin{align*}
\lambda_1&=1/2\,(\mu_{{1}}+\mu_{{4}}+\sqrt {{\mu_{{4}}}^{2}-2\,\mu_{{4}}
\mu_{{1}}+{\mu_{{1}}}^{2}+4\,\mu_{{3}}\mu_{{2}}+4\,\mu_{{3}}})\\
\lambda_2&=\mu_2\\
\lambda_3&=1/2\,(\mu_{{1}}+\mu_{{4}}-\sqrt {{\mu_{{4}}}^{2}-2\,\mu_{{4}}
\mu_{{1}}+{\mu_{{1}}}^{2}+4\,\mu_{{3}}\mu_{{2}}+4\,\mu_{{3}}}).
\end{align*}
Then $Q(\mu)=G\inv A(\ph(\mu))G$, where
$$G=\left[\begin{matrix}1&0\\-1/2\,{\frac {-\mu_{{1}}+\mu_{{4}}+\sqrt {{\mu_{{4}}}^{2}-2\,\mu_{{4}}
\mu_{{1}}+{\mu_{{1}}}^{2}+4\,\mu_{{3}}\mu_{{2}}+4\,\mu_{{3}}}}{1+\mu_{
{2}}}}
&0\end{matrix}\right].$$
This mapping fails to be analytic at zero.  As a consequence, it fails one of the conditions for versality.

One of the nice features of this analytic version of deformation theory is that versal deformations can be characterized in terms of transversality.
Essentially, the idea is that the action of the group at a point
determines a subspace of the tangent space at the point. A deformation
$A(\lambda)$ is transversal if in a neighborhood of $A(0)$, the subspace of the tangent space at each point of $A(\lambda)$ parallel to the family is complementary to the subspace determined by the group action.

One of the nice results in \cite{arnold1} is the determination of the number of parameters of a miniversal deformation of a matrix $A_0$ in terms of the
Jordan decomposition of the matrix.

\begin{thm}[Arnold] If $A_0$ is a matrix, then the number $n$ of parameters of a miniversal deformation of $A_0$ is given by
\begin{equation*}
n=\sum_\lambda n_1+3n_2+5n_3+\cdots,
\label{arneq}
\end{equation*}
where the sum is taken over all eigenvalues $\lambda$ of the matrix, and $n_1\ge n_2\ge\cdots$ are the sizes of the Jordan blocks corresponding to $\lambda$.
\end{thm}
 The basic idea of the proof is as follows. If we consider a one parameter subgroup
$\exp(tB)$ of $\GL(n)$, then a tangent vector to the action of this subgroup is $\frac d{dt}\,\exp(tB)^*(A)=[A,B]$. Consider the map $f$ defined on the $n\times n$ matrices given by $f(B)=[A,B]$. Then the dimension of the tangent space to the group action at $A$ is $\dim(\im(f))$, and so its codimension $\dim(\ker(f))$ is just the dimension of the centralizer of $A$.  This means the dimension of the miniversal deformation is given by the dimension of the centralizer of $A$. Arnold goes on to compute this dimension by looking at a Jordan normal form for $A$, and explicitly constructing its centralizer, from which he concludes the dimension formula above.
\medskip

One should note that the centralizer subspace is not, in general, transverse to the tangent space of the group action, so that one cannot use the centralizer to compute the miniversal deformation directly. Nevertheless, Arnold gives an explicit form for a miniversal deformation based on the Jordan normal form of the matrix.  

A natural question is what is the relation to this idea of transversality and miniversal deformations to the algebraic notion.  First, note that the
cocycle condition $[d,\ph]=0$ is really a condition that at least on the infinitesimal level, the direction indicated by $\ph$ points along the variety of
algebras. This condition does not arise in Arnold's formulation because he was not studying a subvariety of the vector space, so that all tangent directions
are allowed.  Next, the condition
$\ph=[d,\lambda]$ for a coboundary is precisely the same condition as in the Arnold formulation, because it means that $\ph$ is a tangent vector in the
direction of the group action.  Thus $H^2(d)$ represents the transverse directions to the group action.

The above description is more morally correct than precisely correct.  For example, there is a finite dimensional Lie algebra $d$ for which $H^2(d)$ does not
vanish, but there are no nontrivial deformations of $d$. There is an infinitesimal deformation which does not extend to a formal deformation.  Thus, in reality,
there are no curves in the variety of Lie algebra structures through the algebra $d$ which are transverse to the group action.  This means the infinitesimal
picture only gives an approximation to the actual deformation picture.

\section{A stratification of the Moduli space}
Arnold gives a stratification of the space of $n\times n$ matrices in terms of ``types'' of Jordan decompositions.  The authors have been studying moduli spaces of algebras and have found that for small dimensional spaces, the moduli spaces of Lie and associative algebras have a natural stratification by orbifolds of
a very simple type.  For Lie algebras, a portion of the moduli space is given by the action of a group on the space of $n\times n$ matrices, which is not quite
the same action as conjugation by matrices, but is nearly so.  We have an action of $\GL(n)\times\C^*$ on matrices given by
$(G,x)^*(A)=xG\inv A G$. Thus our action is conjugation up to multiplication by a nonzero scalar. When we looked at the picture developed by Arnold, we realized
that a similar stratification of the space of matrices up to conjugation can be given, which is consistent with Arnold's computation of the dimension of the miniversal deformations.  We will describe the decompositions for dimensions 2, 3 and 4 below.
\subsection{Decomposition of $\gl(2)$} The space $\gl(2)$, under the action of $\GL(2)\times\C^*$ has a stratification by the matrix types
\begin{equation*}
A(p:q)=\left[\begin{matrix}p&1\\0&q\end{matrix}\right],\qquad B_1=\left[\begin{matrix}1&0\\0&1\end{matrix}\right],\qquad B_0=\left[\begin{matrix}0&0\\0&0\end{matrix}\right].
\end{equation*}
The first stratum can be parameterized by $\C\P^1/\Sigma_2$, where the action of the symmetric group $\Sigma_2$ on $\C\P^1$ is given by 
permutation of the projective coordinates $(p:q)$. 
In this formulation, the generic point $B(0:0)$ is considered to be an element of $\C\P^1$, except that the number of parameters of the miniversal deformation
of the generic point can be higher.  
The other two strata consists of singleton points.  For the strata $A(p:q)$, the miniversal deformation has one parameter, while the miniversal deformation
of $A(0:0)$ requires 2 parameters. The stratum with the singleton $B_1$ has 3 parameters in the miniversal deformation. Finally, the zero matrix $B_0$ also requires
three parameters for its miniversal deformation.

The space $\gl(2)$, under the action of $\GL(2)$ by conjugation has a stratification
\begin{equation*}
A(p,q)=\left[\begin{matrix}p&1\\0&q\end{matrix}\right],\qquad B(p)=\left[\begin{matrix}p&0\\0&p\end{matrix}\right].
\end{equation*}
The stratum $A(p,q)$ is parameterized by $\C^2/\Sigma_2$, where $\Sigma_2$ acts by permutation of the coordinates $(p,q)$. The stratum $B(p)$ is parameterized by
$\C$. The miniversal deformation of the stratum $A(p,q)$ requires exactly 2 parameters, while the miniversal deformation of the stratum $B(p)$ requires 
4 parameters. 

Let us show how the count of the parameters in the Arnold case is obtained. For an element $A(p,q)$ when $p\ne q$, the parameter $n$ is given by
$n=1+1$, because each eigenvalue contributes 1 to the total. But when $p=q$, we get a Jordan block of size 2, so that $n=2$ again.  For the matrix $B(p)$,
we simply have $n=1+3$ since there are 2 Jordan blocks of size 1.  

The first thing to note is that in our formulation, the matrices $A(p,q)$ can have 2 different Jordan decompositions, so two different strata of Arnold's decomposition are combined to make one complete orbifold stratum.  Secondly, we see that $A(p,q)$ is related to $A(p:q)$ and that the number of parameters
in the projective formulation drops by 1, except for the generic value. A similar relation to $B(p)$ and $B_1$ occurs, since $B_1$ requires 1 less parameter
than $B(p)$ for a miniversal deformation.

\subsection{Decomposition of $\gl(3)$} The space $\gl(3)$, under the action of $\GL(3)\times\C^*$ has a stratification by the matrix types
\begin{equation*}
A(p:q:r)=\left[\begin{matrix}p&1&0\\0&q&1\\0&0&r\end{matrix}\right] B(p:q)=\left[\begin{matrix}p&0&0\\0&p&1\\0&0&q\end{matrix}\right] C_1=\left[\begin{matrix}1&0&0\\0&1&0\\0&0&1\end{matrix}\right]C_0=\left[\begin{matrix}0&0&0\\0&0&0\\0&0&0\end{matrix}\right]
\end{equation*}
The first stratum $A(p:q:r)$ is parameterized by $\C\P^2/\Sigma_3$, where $\Sigma_3$ acts by permuting the projective coordinates $(p:q:r)$. The miniversal deformation
of an element in this stratum has 2 parameters except that $A(0:0:0)$ requires 3 parameters.
The second stratum $B(p:q)$ is parameterized by $\C\P^1$, and there is no action of $\Sigma_2$ on this stratum, which is clear from the fact that $p$ and $q$ play different roles in the matrix.  The miniversal deformation has 4 parameters, except that $B(0:0)$ requires 5 parameters.  The strata $C_1$ and $C_0$ each
require 8 parameters for the miniversal deformation.

The space $\gl(3)$, under the action of $\GL(3)$ by conjugation has a stratification
\begin{equation*}
A(p,q,r)=\left[\begin{matrix}p&1&0\\0&q&1\\0&0&r\end{matrix}\right],\quad B(p,q)=\left[\begin{matrix}p&0&0\\0&p&1\\0&0&q\end{matrix}\right],\quad C(p)=\left[\begin{matrix}p&0&0\\0&p&0\\0&0&p\end{matrix}\right].
\end{equation*}
The first stratum, $A(p,q,r)$ is parameterized by $\C^3/\Sigma_3$, and has a miniversal deformation with 3 parameters.
The second stratum, $B(p,q)$ is parameterized by $\C^2$, and requires 5 parameters for a miniversal deformation.  Finally, the third stratum
$C(p)$ is parameterized by $\C$, and requires 9 parameters.

As in the 2-dimensional example, the number of parameters in the Arnold correspondence is 1 larger than the projective correspondence, except for the generic
element, for which the numbers are the same. Note for example, that the stratum $A(p,q,r)$ consists of the matrices with 3 Jordan blocks and 3 distinct eigenvalues, 2 Jordan blocks, one of size 2 and one of size 1, with 2 distinct eigenvalues, and one Jordan block of size 3, with one eigenvalue. In all three
cases, the formula \refeq{arneq} gives the same number of parameters for the miniversal deformation.

\subsection{Decomposition of $\gl(4)$} The space $\gl(4)$, under the action of $\GL(4)\times\C^*$ has a stratification by the matrix types
\begin{align*}
A(p:q:r:s)=\left[\begin{matrix}p&1&0&0\\0&q&1&0\\0&0&r&1\\0&0&0&s\end{matrix}\right],\quad
E_1=\left[\begin{matrix}1&0&0&0\\0&1&0&0\\0&0&1&0\\0&0&0&1\end{matrix}\right],
E_2=\left[\begin{matrix}0&0&0&0\\0&0&0&0\\0&0&0&0\\0&0&0&0\end{matrix}\right]\\
B(p:q:r)=\left[\begin{matrix}p&0&0&0\\0&p&1&0\\0&0&q&1\\0&0&0&r\end{matrix}\right],
C(p:q)=\left[\begin{matrix}p&1&0&0\\0&q&0&0\\0&0&p&1\\0&0&0&q\end{matrix}\right],
D(p:q)=\left[\begin{matrix}p&0&0&0\\0&p&0&0\\0&0&p&1\\0&0&0&q\end{matrix}\right].
\end{align*}
The first stratum $A(p:q:r:s)$ is parameterized by $\C\P^3/\Sigma_4$, and has a miniversal deformation with 3 parameters, except for $A(0:0:0:0)$ which 
requires 4.  The second stratum $B(p:q:r)$ is parameterized by $\C\P^2/\Sigma_2$, where $\Sigma_2$ acts by interchanging the coordinates $q$ and $r$. This
stratum requires 5 parameters, except for $B(0:0:0)$ where we need 6.  The third stratum $C(p:q)$ is parameterized by $\C\P^1/\Sigma_2$, and requires 7 parameters, except for $C(0:0)$ which requires 8.  The fourth stratum $D(p:q)$ is parameterized by $\C\P^1$, and requires 9 parameters, except for $D(0:0)$,
which requires 10.  Finally, the strata $E_1$ and $E_0$ require 15 parameters.

The space $\gl(4)$, under the action of $\GL(4)$ has a stratification by the matrix types
\begin{align*}
A(p,q,r,s)=\left[\begin{matrix}p&1&0&0\\0&q&1&0\\0&0&r&1\\0&0&0&s\end{matrix}\right],\quad
E_(p)=\left[\begin{matrix}p&0&0&0\\0&p&0&0\\0&0&p&0\\0&0&0&p\end{matrix}\right]\\
B(p,q,r)=\left[\begin{matrix}p&0&0&0\\0&p&1&0\\0&0&q&1\\0&0&0&r\end{matrix}\right],
C(p,q)=\left[\begin{matrix}p&1&0&0\\0&q&0&0\\0&0&p&1\\0&0&0&q\end{matrix}\right],
D(p,q)=\left[\begin{matrix}p&0&0&0\\0&p&0&0\\0&0&p&1\\0&0&0&q\end{matrix}\right].
\end{align*}
The stratification is as follows: $A(p,q,r,s)$ is parameterized by $\C^4/\Sigma_4$, $B(p,q,r)$ is parameterized by $\C^3/\Sigma_2$, $C(p,q)$ is parameterized by
$\C^2/\Sigma_2$, $D(p,q)$ is parameterized by $\C^2$, and finally $E(p)$ is parameterized by $\C$. In each case, the number of parameters is 1 larger than the generic number in the projective picture.
\subsection{How to understand the strata}
As the examples have shown, it is always possible to combine several of the strata given by Jordan types in the Arnold model to give a single stratum which has a nice description as an orbifold. What is not so obvious is that the deformation picture is completely captured by the orbifold stratification, by which we mean
that deformations either occur along a stratum, are given by a jump
deformation to a point in another stratum, or ``factor through a jump
deformation'' to a point in another stratum.  By ``factoring through a
jump deformation" we mean that a deformation goes along a neighborhood of a point in the stratum to which there is a jump deformation. A jump deformation from $A$ to $B$ is a deformation $A(\lambda)$ for which $A(0)=A$ and $A(\lambda)$ is equivalent to $B$ for $\lambda\ne 0$. 

A simple example of a jump deformation can be given for $2\times 2$ matrices. The deformation $A(t)=\left[\begin{smallmatrix}1&t\\0&1\end{smallmatrix}\right]$ is a jump deformation from $A=\left[\begin{smallmatrix}1&0\\0&1\end{smallmatrix}\right]$ to $B=\left[\begin{smallmatrix}1&1\\0&1\end{smallmatrix}\right]$, because $A(t)$ is equivalent to $B$ for all $t\ne0$.

We also noticed the following.  When $n$ is even, all the strata have miniversal deformations with an even number of parameters, and when $n$ is odd, they all
have an odd number of parameters. This fact is easy to prove. In our examples, the number of parameters of a miniversal deformation were different for each stratum. However, this is an artifact of the size of our matrices, and is not true in general.  For $n=6$, the following two matrices represent different strata, but each of them has a miniversal deformation with 12 parameters.
\begin{equation*}
\left[\begin{matrix}p&0&0&0&0&0\\0&p&0&0&0&0\\0&0&p&1&0&0\\0&0&0&q&1&0\\0&0&0&0&r&1\\0&0&0&0&0&s\end{matrix}\right],\qquad
\left[\begin{matrix}p&0&0&0&0&0\\0&q&0&0&0&0\\0&0&r&1&0&0\\0&0&0&p&0&0\\0&0&0&0&q&0\\0&0&0&0&0&r\end{matrix}\right].
\end{equation*}

The construction of the strata can be accomplished in a straightforward manner, with one stratum for each partition of the number $n$. Thus, the number of strata of the moduli space of $n\times n$ matrices under the action of $\GL(n)$ by conjugation is exactly the number of partitions of $n$.

It is important to understand that the stratification we have given is not simply a decomposition of the space in terms of pieces which consist of elements whose miniversal deformations have the same number of parameters, but in fact, the deformation theory is consistent with this stratification and no other.


\end{document}